\documentclass[%
 reprint,
 amsmath,amssymb,amsfont
 prl,
]{revtex4-2}

\usepackage[utf8]{inputenc}
\usepackage[english]{babel}
\usepackage{graphicx}
\usepackage{dcolumn}
\usepackage{bm}
\usepackage{comment}
\usepackage[table,dvipsnames]{xcolor}
\usepackage{tikz}
\usepackage{subcaption}
\usepackage{siunitx}
\usepackage{booktabs}
\usepackage{standalone}
\usepackage{xr-hyper} 
\usepackage{hyperref}
\usepackage{xcolor}
\hypersetup{%
    colorlinks,
    linkcolor={blue!50!black},
    citecolor={blue!50!black},
    urlcolor={blue!50!black}
}

\makeatletter
\newcommand*{\addFileDependency}[1]{
  \typeout{(#1)}
  \@addtofilelist{#1}
  \IfFileExists{#1}{}{\typeout{No file #1.}}
}
\makeatother

\newcommand*{\myexternaldocument}[1]{%
    \externaldocument{#1}%
    \addFileDependency{#1.tex}%
    \addFileDependency{#1.aux}%
}
\myexternaldocument{appendix}

\usepackage{multirow}

\newcommand{\HP}{\mathbf{P}} 
\newcommand{\SF}{\mathcal{F}} 
\newcommand{\bF}{\mathbf{F}}
\newcommand{\uu}{\mathbf{u}}
\newcommand{\UU}{\mathbf{U}}
\newcommand{\tuu}{\tilde{\mathbf{u}}}
\newcommand{\NN}{\vec{N}}
\newcommand{\XX}{\vec{X}} 
\newcommand{\dd}{\partial}
\newcommand{\HH}{\vec{H}}
\newcommand{\eps}{\epsilon}
\newcommand{\reals}{\mathbb{R}}
\DeclareMathOperator{\Div}{div}
\DeclareMathOperator{\tr}{tr}
\DeclareMathOperator{\diag}{diag}

\begin{document}

\UseRawInputEncoding 

\preprint{APS/123-QED}

\title{Scroll Waves and Filaments in excitable Media of higher spatial Dimension}

\author{Marie Cloet\textsuperscript{1,2}}
\author{Louise Arno\textsuperscript{1,2}}%
\author{Desmond Kabus\textsuperscript{1,3,2}}
\author{Joeri Van der Veken\textsuperscript{4}}
\author{Alexander V. Panfilov\textsuperscript{5,6,7}}
\author{Hans Dierckx\textsuperscript{1,2}}
 \email{h.dierckx@kuleuven.be}
\affiliation{%
 \textsuperscript{1}Department of Mathematics, KU Leuven Campus Kortrijk (KULAK), Kortrijk 8500, Belgium, 
 \textsuperscript{2}iSi Health, Institute of Physics-based Modeling for In Silico Health, KU Leuven, Leuven 3000, Belgium,
 \textsuperscript{3}Laboratory of Experimental Cardiology, Leiden University Medical Center (LUMC), Leiden 2333 ZA, The Netherlands, 
 \textsuperscript{4}Department of Mathematics, KU Leuven, Leuven 3001, Belgium, 
 \textsuperscript{5}Department of Physics and Astronomy, Ghent University, Ghent 9000, Belgium, 
 \textsuperscript{6}Laboratory of Computational Biology and Medicine, Ural Federal University, Ekaterinburg 620002, Russia, 
 \textsuperscript{7}World-Class Research Center ``Digital biodesign and personalized healthcare'', Sechenov University, Moscow 119991, Russia
}%




\date{Received 2 May 2023; accepted 30 August 2023}

\begin{abstract}
Excitable media are ubiquitous in nature, and in such systems the local excitation tends to self-organize in traveling waves, or in rotating spiral-shaped patterns in two or three spatial dimensions. Examples include waves during a pandemic or electrical scroll waves in the heart. 
Here we show that such phenomena can be extended to a space of four or more dimensions and 
propose that connections of excitable elements in a network setting can be regarded as additional spatial dimensions. Numerical simulations are performed in four dimensions using the FitzHugh-Nagumo model, showing that the vortices rotate around a two-dimensional surface which we define as the superfilament.
Evolution equations are derived for general superfilaments of co-dimension two in an $N$-dimensional space and their equilibrium configurations are proven to be minimal surfaces. We suggest that biological excitable systems, such as the heart or brain which have non-local connections can be regarded, at least partially, as multidimensional excitable media and discuss further possible studies in this direction.
\end{abstract}

\maketitle


\textit{Introduction.} Many real-world systems exhibit a large non-linear response to external or internal stimuli, and are therefore called excitable systems. In the continuum limit, several of these systems are reasonably understood, as unifying geometric principles have been revealed that enable us to understand and quantify their dynamics. As a first example, non-linear traveling waves of constant amplitude have been observed as depolarization waves in neural~\cite{bressloff2014} and cardiac tissue~\cite{clayton2011}, combustion~\cite{zeldovich1938} or oxidation~\cite{rovinskii1984} processes, Mexican waves in crowds~\cite{farkas2002}, and infection waves in a pandemic~\cite{medlock2003}. This extremely diverse range of phenomena has been described mathematically in the continuum limit using reaction-diffusion equations
\begin{equation}\label{RDE}
    \partial_t\uu= \HP \Delta_N \uu + \bF(\uu)\,,
\end{equation}
where $\uu(\vec{r},t)$ is a state vector with $S$ components that varies in $N$-dimensional space and time, and $\HP$ is a diagonal $S \times S$ matrix containing the diffusivity of each state variable, $\Delta_N$ the Laplacian operator in $N$ spatial dimensions and $\bF(\uu)$ a local excitation model. The state vector can support waves and it was shown mathematically that their wave fronts undergo curvature-driven dynamics, with the lowest order effect behaving like surface
tension~\cite{kuramoto1980,keener1988,dierckx2011}. These curvature effects are governed by the medium characteristics, mathematically represented by the model parameters. Wave fronts are known to be lines in two-dimensional (2D) media, and surfaces in 3D media. This property can be summarized as wave fronts having spatial codimension $C=1$, i.e., they are $(N-C)$\nobreakdash-dimensional structures in an $N$-D space. 

In observations of excitable surfaces ($N=2$), excitation patterns were often found to organize into rotating spiral-shaped patterns. Their study is largely motivated by cardiological applications, since electrical rotors have been observed on the surface of the heart during rhythm disorders. Like for wave fronts, spiral waves drift is affected by geometric factors, such as the local Gaussian curvature of the excitable surface~\cite{mikhailov1994}, possibly incorporating anisotropy of wave propagation within the surface~\cite{dierckx2013}. Spiral waves on a surface rotate around a core region, which is commonly idealized into a phase singularity (PS) point, of co-dimension 2~\cite{gray1998}. The case where quasi-periodically moving (meandering) PSs occur is more involved, and has been described recently~\cite{tomii2021,arno2021,kabus2022}. 
In 3D media, such as cardiac tissue, a stack of spiral waves is known as a scroll wave. A scroll wave rotates around a curve consisting of PSs, the filament which has dimension $N-C = 3-2 = 1$. The curvature-induced dynamics of filaments was first established numerically in~\cite{panfilov-rudenko1987,panfilov1984,panfilov1987}
for isotropic excitable media 
where they were shown to exhibit so-called filament tension~\cite{biktashev1994}, a medium-dependent emerging parameter. For positive filament tension, filaments straighten up~\cite{panfilov1984,panfilov1987}, while for negative filament tension, filaments elongate and tend to break up~\cite{panfilov-rudenko1987}, if the medium is thick enough~\cite{dierckx2012}. In the context of cardiac arrhythmias, the topological configuration of filament curves in the cardiac wall (e.g.,
$O$, $I$ or $U$-shaped) is being used as a classifier for arrhythmia patterns. Moreover, spatio-temporal drift (due to internal interaction or gradients in the medium) is thought to be responsible for transitions between arrhythmias~\cite{jalife1998}. 

From the above overview on the current knowledge, a natural question arises that is at first sight purely mathematical: Which other patterns of co-dimension $C$ exist in $N$-dimensional media? In this paper, we provide some answers to these questions. We confirm that indeed other patterns exist, more specifically, superfilaments ($C=2$), i.e., the organizing center of superscrolls in a space with $N>3$ spatial dimensions. 

\begin{figure}
    \centering
        \begin{subfigure}[t]{0.15\linewidth}
            \centering
            \begin{tikzpicture}[scale=0.6]
                \draw[dashed] (0,-1.5) -- (0,1.5);
                \draw[thick] (0,-1) -- (0,1);

                \fill[black, draw=black, thick] (0,0) circle (2pt);
                \fill[black, draw=black, thick] (0,1) circle (2pt);
                \fill[black, draw=black, thick] (0,-1) circle (2pt);
            \end{tikzpicture}
            \caption{1D}
            \label{fig: 1D stencil}
        \end{subfigure}
        \hfill
        \begin{subfigure}[t]{0.25\linewidth}
            \centering
            \begin{tikzpicture}[scale=0.6]
                \draw[dashed] (-1.5,0) -- (1.5,0);
                \draw[dashed] (-1.5,-1) -- (1.5,-1);
                \draw[dashed] (-1.5,1) -- (1.5,1);
                \draw[dashed] (0,-1.5) -- (0,1.5);
                \draw[dashed] (-1,-1.5) -- (-1,1.5);
                \draw[dashed] (1,-1.5) -- (1,1.5);
                \draw[thick] (-1,0) -- (1,0);
                \draw[thick] (0,-1) -- (0,1);
    
                \fill[black, draw=black, thick] (-1,-1) circle (1pt);
                \fill[black, draw=black, thick] (-1,1) circle (1pt);
                \fill[black, draw=black, thick] (1,-1) circle (1pt);
                \fill[black, draw=black, thick] (1,1) circle (1pt);
                \fill[black, draw=black, thick] (-1,0) circle (2pt);
                \fill[black, draw=black, thick] (0,0) circle (2pt);
                \fill[black, draw=black, thick] (1,0) circle (2pt);
                \fill[black, draw=black, thick] (0,1) circle (2pt);
                \fill[black, draw=black, thick] (0,-1) circle (2pt);
            \end{tikzpicture}
            \caption{2D}
            \label{fig: 2D stencil}
        \end{subfigure}
        \begin{subfigure}[t]{0.45\linewidth}
            \centering
            \begin{tikzpicture}[scale=0.8]
                \draw[dashed] (-1,0) -- (1,0);
                \draw[dashed] (-1,-1) -- (1,-1);
                \draw[dashed] (-1,1) -- (1,1);
                \draw[dashed] (0,-1) -- (0,1);
                \draw[dashed] (-1,-1) -- (-1,1);
                \draw[dashed] (1,-1) -- (1,1);
                \draw[thick] (-1,0) -- (1,0);
                \draw[thick] (0,-1) -- (0,1);
                \draw[thick] (-0.6,-0.4) -- (0.6,0.4);
             
                \foreach \y in {-1,0,1} {
                     \draw[dashed] (-.6,-0.4+\y) -- (0.6,0.4+\y);
                     }
                \foreach \x in {-1,1} {
                     \draw[dashed] (-.6+\x,-0.4) -- (0.6+\x,0.4);
                     }
                 
                \draw[dashed] (-0.6,-1.4) -- (-0.6,0.6);
                \draw[dashed] (0.6,-0.6) -- (0.6,1.4);
                \draw[dashed] (-1.4,-0.4) -- (0.4,-0.4);
                \draw[dashed] (-0.4,0.4) -- (1.6,0.4);
    
                \fill[black, draw=black, thick] (-1,-1) circle (1pt);
                \fill[black, draw=black, thick] (-1,1) circle (1pt);
                \fill[black, draw=black, thick] (1,-1) circle (1pt);
                \fill[black, draw=black, thick] (1,1) circle (1pt);
                \fill[black, draw=black, thick] (-1.6,-0.4) circle (1pt);
                \fill[black, draw=black, thick] (0.4,-0.4) circle (1pt);
                \fill[black, draw=black, thick] (-.6,-1.4) circle (1pt);
                \fill[black, draw=black, thick] (.6,-0.6) circle (1pt);
                \fill[black, draw=black, thick] (-.6,0.6) circle (1pt);
                \fill[black, draw=black, thick] (0.6,1.4) circle (1pt);
                \fill[black, draw=black, thick] (-0.4,.4) circle (1pt);
                \fill[black, draw=black, thick] (1.6,.4) circle (1pt);
                \fill[black, draw=black, thick] (-1,0) circle (2pt);
                \fill[black, draw=black, thick] (0,0) circle (2pt);
                \fill[black, draw=black, thick] (1,0) circle (2pt);
                \fill[black, draw=black, thick] (0,1) circle (2pt);
                \fill[black, draw=black, thick] (0,-1) circle (2pt);
                \fill[black, draw=black, thick] (-0.6,-0.4) circle (2pt);
                \fill[black, draw=black, thick] (0.6,0.4) circle (2pt);
                
                \draw[-latex] (-2.5,1) -- (-2,1) node[right] {$y$};
                \draw[-latex] (-2.5,1) -- (-2.5,1.5) node[right] {$z$};
                \draw[-latex] (-2.5,1) -- (-2.8,0.8) node[below] {$x$};
            \end{tikzpicture}
            \caption{3D}
            \label{fig: 3D stencil}
        \end{subfigure}
    \caption{
     Stencil for the second order accurate discretization of the Laplacian. Nodes connected to the central point obtain a weight $1/h^2$, see Eq.~\eqref{RDnetwork}.
    }
    \label{fig: Laplacian stencil}
\end{figure}
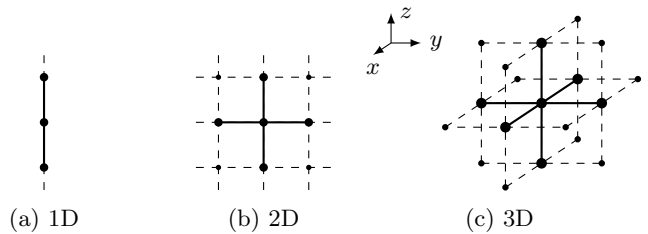

\textit{Research hypothesis. } We claim that such higher dimensional patterns are significant for real-world excitable systems. The crucial point is that simulations on an $N$\nobreakdash-D grid are related to a regular network of connected excitable elements. Indeed, when the aforementioned systems in $N=2$ or $N=3$ spatial dimensions are studied numerically, they are often simulated on a Cartesian grid with lattice constant $h$. Using the second order accurate Laplacian stencil (see Fig.~\ref{fig: Laplacian stencil}), the discretization of Eq. (1) for the $i$-th node in the grid becomes
\begin{equation}
    \partial_t\uu_i = \HP \sum_{j : A_{ij}=1} \frac{\uu_j - \uu_i}{h^2}   + \bF(\uu_i)\,. \label{RDnetwork}
\end{equation}
Here, $A_{ij} = 1$ if the node $i$ is next to $j$ and $0$ otherwise. In  $\reals^1$, every interior point will have $2$ such neighbors, in  $\reals^2$,  $4$ neighbors, etc. Thus, if one will discretize similarly the Laplacian in  $\reals^4$, every interior point will have $8$ neighbors, and in $\reals^N$, one obtains $2N$ such neighbors. Note, however, that Eq.~\eqref{RDnetwork} can also be regarded as an excitable network~\cite{nakao2010}, where $A_{ij}$ is the adjacency matrix of the graph consisting of the edges that form the network. In this view, going to $N > 3$ actually makes sense, as an idealized model for excitable networks in which each node is connected to $2N$ neighbors. In network language, such graph would be a limit of a translationally invariant network with average degree $2N$, also called a regular network. In view of real-world networks with unexplained emergent pattern formation, $N$ can grow large. For example, on social networks, the average degree is of magnitude $10^2-10^3$~\cite{myers2014}. In brain tissue, cerebellar granule cells make about 675 synaptic connections with Purkinje cells, which in their turn receive input from around 185,000 parallel fibers~\cite{harvey1991}. If a network with high degree is realized in 3D space that surrounds us, it implies that the network exhibits long-range, or global connections. 
Within this work, analytical results are obtained for any integer $N>2$. However, we will consider in simulations only one extra dimension, such that $N=4$, which already shows interesting new physical phenomena.  

\textit{Analytical derivation.} Motivated by multiple applications, as presented above, and inspired by previous results in 2 and 3 dimensions \cite{keener1988,biktashev1994,verschelde2007}, the equation of motion for the superfilaments of co-dimension $C=2$ will be derived here. In those works, a scroll wave was constructed by stacking a set of spiral waves along a filament curve. Let $x^i$, $1\leq i\leq N$, be the coordinates in an $N$-dimensional domain. By definition, the superfilament
$\SF$ will be an
$(N-2)$-dimensional submanifold, hence we can define its coordinates $\sigma^m$, where $1\leq m\leq N-2$. The embedding of the superfilament in the $N$-dimensional domain is then given by $X^i(\sigma^1,\ldots,\sigma^{N-2})$, $1\leq i\leq N$.
In one point of $\SF$, we construct 2 orthogonal unit normal vectors $\vec{N}_c$, with $c\in\{1,2\}$. The spatial change of these vectors between different points on the submanifold is given by the Weingarten formula \cite{lee2018}: 
\begin{align}
D_m \vec{N}_c = - S_{\vec{N}_c}\left(\vec{e}_m\right) + \nabla_m^{\perp}\vec{N}_c\,. \label{normalconnection}
\end{align}
Here, $D_m$ is the derivative in $\reals^N$ in the direction of coordinate $\sigma^m$, $S_{\vec{N}_c}$ is the shape operator inherited from the normal direction $\vec{N}_c$ and $\vec{e}_m$ is the tangent vector to $\mathcal{F}$ corresponding to coordinate direction $\sigma^m$. The second term $\nabla_m^{\perp}\vec{N}_c$ is the normal connection, describing the relative rotation of both normals when one moves along the superfilament. This term can be set to 0 by choosing the normal vectors $\vec{N}_1$ and $\vec{N}_2$ to be parallel with respect to the normal connection. This procedure is similar to choosing a relative parallel frame to a curve in 3D \cite{bishop1975}. 

Next, we introduce normal coordinates $\rho^c$ for points in the vicinity of the superfilament:
\begin{align}\label{eq: choice of coordinates}
 x^i = X^i\left(\sigma^m\right) + \sum_{c=1}^2\rho^c N_c^i.
\end{align}
This allows us to define the normal vectors $\vec{N}_c$ also in the vicinity of the superfilament, by $N_c^i = \frac{\partial x^i}{\partial \rho^c}$. The divergence of $\vec{N}_c$ is then related to the curvature of the superfilament. Indeed, by construction, $D_d \vec{N}_c =\vec{0}$ for $c\in\{1,2\}$, $d\in\{N-1,N\}$, such that with Eq.~\eqref{normalconnection} we obtain~\cite{lee2018,dajczer2019}:
\begin{equation}\label{eq: definition divergence}
\Div\left(\vec{N}_c\right)=\sum_{m=1}^{N-2} \vec{e}_m\cdot D_m \vec{N}_c = -\tr\left(S_{\NN_c}\right) = - \vec{N}_c \cdot \vec{H} \,,
\end{equation}
with $\vec{H}$ the (non-normalized) mean curvature of the superfilament.

With these preparations, we can now explicitly construct the superscroll as a stack of spiral wave solutions. Say $\UU(x,y; \phi)$ is a spiral wave solution in the 2D plane, rotated around the origin over an angle $\phi$. Then, we propose as a solution to Eq. \eqref{RDE} in $N$ dimensions:
\begin{align}
     \uu\left(x^i, t\right) &= \UU\left(\rho^1, \rho^2, \phi(\sigma^m, t)\right) + \tuu\left(x^i, t\right), \nonumber \\
    \dot{\phi}\left(\sigma^m, t\right) &= \omega + \tilde{\omega}\left(\sigma^m, t\right)
\end{align}
where $\omega$ is the spiral's natural rotation frequency in a plane. 
This approximation is a gradient expansion around a planar superfilament, 
and the $\tuu$ and $\tilde{\omega}$ terms are higher-order corrections in the extrinsic curvature of the superfilament, 
similar to earlier work on classical scroll wave filaments \cite{keener1988,verschelde2007}.
Substituting the expressions in the reaction diffusion equation (RDE) brings, in leading order,
\begin{equation}\label{eq: linear RDE}
\dot{\tuu} - \hat{\mathbf{L}}
\tuu=
    \sum_{c=1}^2 \left[\HP\Div(\NN_c)
    + \left(\NN_c\cdot \dot{\XX}\right)\right]\frac{\partial\UU}{\partial \rho^c}
    + \omega \frac{\partial \UU}{\partial \phi},
\end{equation}
where $\hat{\mathbf{L}} = \HP \Delta_2 + \omega \dd_\phi + \bF'(\UU)$. This operator only depends on the standard 2D unperturbed rotating spiral solution. From symmetry arguments, it is known to have critical eigenmodes (Goldstone modes)~\cite{keener1988,barkley1994,hakim2002,biktasheva2003,verschelde2007,biktasheva2009}, as well as response functions~\cite{biktasheva2003}. They are critical modes of the adjoint operator $\hat{\mathbf{L}} ^\dagger = \HP^H \Delta_2 -\omega \dd_\phi + \bF'^H(\UU)$. Taking the overlap integral of Eq. \eqref{eq: linear RDE} with the response functions delivers, after averaging over one rotation period~\cite{biktashev1994,verschelde2007}: 
\begin{equation}
\dot{\XX} = -\gamma_1 \sum_{c=1}^2\Div\left(\NN_c\right)\NN_c - \gamma_2 \sum_{c,d=1}^2\varepsilon_{cd}\Div\left(\NN_d\right)\NN_c, 
\end{equation}
Here, $\gamma_1$ and $\gamma_2$ are the scalar and pseudoscalar filament tension coefficients that were introduced in \cite{biktashev1994} and $\varepsilon_{cd}$ are the components of the Levi-Civita symbol $\bar{\bar{\mathbf{\varepsilon}}} = \left(\begin{smallmatrix} 0 & 1\\ -1 & 0\end{smallmatrix} \right)$. Using Eq. \eqref{eq: definition divergence} we can express the superfilament dynamics in terms of its mean curvature:
\begin{equation}
\dot{\XX}= \gamma_1\HH +  \gamma_2 \bar{\bar{\mathbf{\varepsilon}}} \cdot \HH, 
\label{eom_h}
\end{equation}
Our result is consistent with the lower-dimensional case: For $N=3$, the superfilament becomes a filament curve with arc length $s$, curvature $k$, tangent vector $\vec{T}$ and unit normal $\vec{N}$, such that $\vec{H} = k \vec{N} = \frac{\dd \vec{T}}{\dd s}$. Then, Eq. \eqref{eom_h} reduces to the classical result of Biktashev et al. \cite{biktashev1994}. 

\begin{figure*}
    \begin{tabular}{cc}
    \raisebox{2cm}{a)}
    \includegraphics[width=0.6\linewidth]{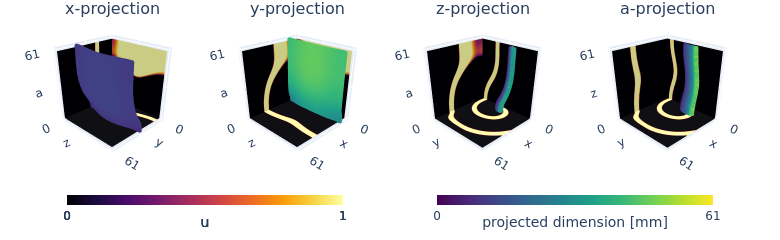}
    & \multirow[t]{2}{*}{
        \raisebox{2cm}{c)}
        \raisebox{-0.5\height}{\includegraphics[trim=210 100 210 100, clip, width=0.25\textwidth]{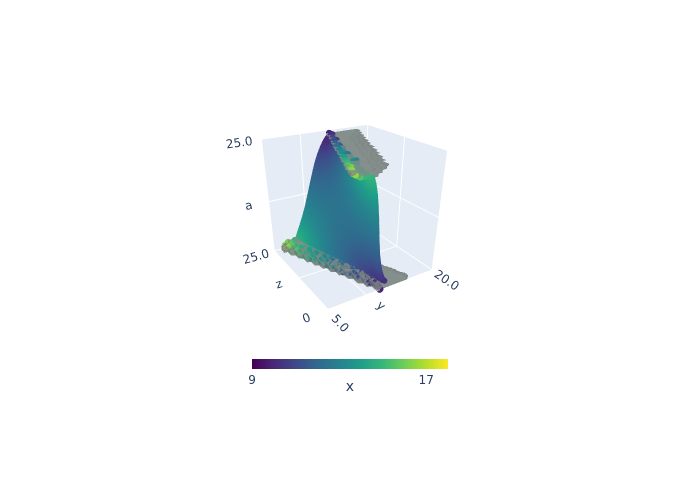}}
        } \\
    \raisebox{2cm}{b)}
    \includegraphics[width=0.6\linewidth]{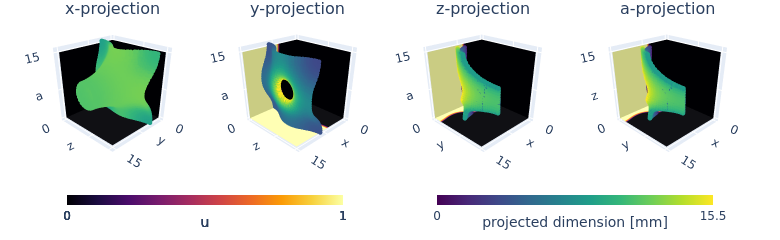}
    &
    \end{tabular}
    \caption{Superscrolls rotating around superfilaments in 4 spatial dimensions. Projections of superfilament points are shown centrally; the back left, back right and bottom of the medium are colored according to the $u$-value at those planes. (a) A superscroll rotating around a 2D curved superfilament ($n=61, h=1\,$mm). Four visualizations of the same timeframe are shown, each frame with a different dimension rendered as color of the superfilament. 
    (b) A pierced superfilament ($n=31, h=0.5\,$mm), shown with projections as in (a). See Supplemental Material~\cite{suppmat} for interactive figures and movies of these simulations. 
    (c) Projection of superfilament points on the Y, Z, and A-axes that spans a minimal surface between obstacles (shaded in grey) at opposite domain boundaries 
    ($n=50, h=0.5\,$mm).
    }
    \label{fig:curvedsuperfilament}
\end{figure*}

From the equation of motion \eqref{eom_h}, it follows that stationary superfilaments must have $\vec{H} = \vec{0} $, thus the total mean curvature must vanish. For the case $N=4$, this means that a stationary superfilament must be a minimal surface, generalizing Wellner's minimal principle to superscrolls  \cite{wellner2002}. It proves a property that we postulated before based on an action principle \cite{dierckx2015}, see also the discussion on Fig. \ref{fig:curvedsuperfilament}c below. 

Similar to wave fronts and classical filaments, the motion of superfilaments in a homogeneous medium is driven by curvature. It is possible to explicitly show that in the case of positive filament tension, the superfilament will decrease its size, as follows. 
With $A$ the total area of a (hyper)surface $\mathcal{S}$ and $\dot{\vec{X}}$ prescribing temporal evolution for all points on that surface, Gauss' first law for the variation of surface area states 
$dA/dt=
- \int_\SF \HH \cdot \dot{\vec{X}} \mbox{d}V$,
with $\mbox{d}V$ the volume form on the (hyper)surface \cite{frankel1997}. 
For superfilaments $\mathcal{F}$, combining this relation with Eq. \eqref{eom_h} yields:
\begin{align}\label{dSdt}
    \frac{dA}{dt}
    & = - \int_\SF \HH \cdot \left[ \gamma_1 \HH + \gamma_2 \bar{\bar{\mathbf{\varepsilon}}} \cdot \HH \right] \mbox{d}V
    = - \gamma_1 \int_\SF \lVert \vec{H} \rVert^2 \mbox{d}V
\end{align}
with $\lVert\,.\,\rVert$ the Euclidean norm of a vector.
Eq.~\eqref{dSdt} gives the analogue of the change-of-length rate $dL/dt = - \gamma_1 \int k^2 ds$ for classical filaments~\cite{biktashev1994}, which we retrieve as the case $N=3$. We have here shown that superfilaments monotonically decrease their generalized surface area in case of positive filament tension. As a corollary, stable superscrolls can persist in excitable media as minimal surfaces (or their higher-dimensional counterparts), e.g., by anchoring to opposite boundaries of the medium, as in~\cite{wellner2002,dierckx2015}. Stable knotted structures may also exist~\cite{keener1992}, but this investigation falls outside our present scope. Conversely, if $\gamma_1 <0$, classical filaments are known to destabilize, which is named as a possible pathway from ventricular tachycardia in the heart to fibrillation~\cite{panfilov1998, dierckx2012}. Based on Eqs.~\eqref{eom_h},\eqref{dSdt}, this phenomenon can also occur in higher-dimensional media.  

\textit{Numerical results.} To substantiate the analytical results, the reaction-diffusion system~\eqref{RDE} 
was integrated numerically using the Forward Euler method, 
with an $N$-dimensional hypercube as the domain: $\vec{r} \in {[0,L]}^N$ and no-flux boundary conditions at the edges. 
In each dimension, space was discretized using $n$ points along each dimension. 
This discretization results in a lattice with spacing $h=\frac{L}n$, chosen small enough such that no discretization artifacts could be observed in the state variable
fields. 
For the local excitation dynamics  $\bF(\uu)$, we
chose the FitzHugh-Nagumo (FHN) kinetics~\cite{fitzhugh1961}: $\uu = {[u,v]}^T$, $\bF ={[f(u,v), g(u,v)]}^T$ with 
\begin{align}
f&= \eps^{-1}\left(u+\frac{u^3}{3}-v\right), & g&=\eps\left(u-av+b\right)
\end{align} 
and $\HP = \diag(1,0)$, meaning that only $u$ is diffused. We used $\eps=0.3$, $a=0.5$, and $b=0.68$, which corresponds to the $\gamma_1>0$ regime in 3D. 
The FHN model was originally formulated for neural conduction, but has been applied equally to continuous models of cardiac excitation~\cite{biktasheva2006}. As such, it is well suited to explore the boundary between discrete and continuous excitation models. 

We here report numerical results for $N=4$ only. In a medium of $61 \times 61 \times61 \times61$ voxels of size $h=1\,$mm, a traveling wave was initiated by setting $u=1$ inside a ball around the point $(-\SI{20}{\milli\meter},-\SI{20}{\milli\meter},\SI{31}{\milli\meter},\SI{45}{\milli\meter})$ of radius $\SI{40}{\milli\meter}$. After $\SI{15}{\milli\meter}$, $u=1$ was set for $y\leq \SI{10}{\milli\meter}$. This protocol resulted in a curved planar superfilament, to which we applied filament
tracking in each pair of coordinate planes \cite{fenton1998} with threshold values $u=0$, $v=0$. We thus created the first superscroll
rotating around a superfilament surface, see Fig.~\ref{fig:curvedsuperfilament}a. We show only 3D projections here. Other snapshots and movies are given in the Supplemental Material~\cite{suppmat}. 



Secondly, we created a non-trivial state with a pierced superfilament, see panel~\ref{fig:curvedsuperfilament}b. We started from a stable superscroll with planar superfilament in ${[\SI{0}{\milli\meter},\SI{61}{\milli\meter}]}^4$, limited the domain to ${[\SI{0}{\milli\meter}, \SI{15}{\milli\meter}]}^4$, doubled spatial resolution to $h=\SI{0.5}{\milli\meter}$ by linear interpolation, and set $u=0$ in the hypercyilinder ${(z-8)}^2 + {(a-8)}^2 < 25$.  
Remarkably, the shrinking of the hole in the filament does not occur in the plane of the superfilament, but according to Eq. \eqref{eom_h} at an angle to it, since the pseudoscalar filament tension component $\gamma_2 \neq 0$ for this system, reminiscent of the oblique drift of classical filaments under an external perturbation \cite{biktashev1994, dierckx2009}.


In a third simulation, we created a 
superfilament in 4D that is anchored to 2 obstacles at opposite domain boundaries, which shows indeed that a minimal surface spanned between the obstacles is a stationary configuration. See detailed set-up in the Supplemental Material~\cite{suppmat}.

\textit{Discussion.}
We have shown via both analytical calculations and simulations that adding a fourth spatial dimension to a reaction-diffusion systems yields new structures: a superscroll and a superfilament. 
It is expected that even more complex structures will arise when the number of dimensions $N$, or the co-dimension $C$ of the structure is further increased. 
Our motivation for this research is the link between Laplacian stencils used for simulations according to the Euler method (see Fig.~\ref{fig: Laplacian stencil}) and regular networks~\cite{nakao2010}. 
When the classical scroll waves and filament were introduced in cardiology, they enabled a better understanding of the emergent patterns seen at the outer surface of the heart.
Also, the distance-weighed averaging of the excitation patterns is reflected in the electrical recordings near the heart or on the body surface, and the properties of the scroll wave are now known to affect the frequencies present in the electrograms~\cite{qu2004}.  

With this Letter, we aim to open up research on scrolls in more dimensions, here called superscrolls, in relation to non-linear excitation dynamics in networks. Several additional facts will have to be accounted for. First, we assume here a regular, translationally invariant network. In realistic networks, this condition is typically not met. However, locally there will be closed loops along which the excitation dynamics can travel and re-excite; we leave open the possibility that its
description may require also a non-integer number of dimensions~\cite{wen2021}. Secondly, the properties of a network typically vary between the excitable elements. Heterogeneity of excitable nodes can be dealt with as a gradient as in cardiac modeling, whereas the local variations in nodal degree (i.e., number of connections) may require a probabilistic approach. Thirdly, in our derivation of curvature-driven dynamics, a linear coupling between elements was assumed in  Eq.~\eqref{RDnetwork}. Such effective coupling is also traditional in cardiac tissue modeling, although it is known that signal conduction consists of a traveling wave along the cell membrane interleaved with conduction over gap junctions.

Finally, we suggest to search for superscrolls in real-life systems, as they may offer geometric insight in yet unknown mechanisms. In this search, the notion of excitability can be taken very broadly, as any arrangement of elements that responds non-linearly to stimuli. Examples in human society may include spreading of opinions or epidemics. In technology, contamination, combustions, (temporary) failure and overflow may provide examples. Also oscillatory media may host superscrolls, but the distinction with phase waves needs to be made. 

\textit{Conclusion.} In this work, we have introduced a mathematical generalization of scroll waves and filaments to $N$-dimensional Euclidean space. Superfilaments were shown to obey curvature-driven dynamics and an area-minimizing principle in case of positive tension. Although the world around us is only three-dimensional, we conjecture that superscrolls can occur in media where oscillatory or excitable elements are connected on average to more than
$2N$~neighbors, i.e., in networks with non-local connections, such as neural tissues, and social, biochemical or technological networks.

\begin{acknowledgments}
A. V. Panfilov is grateful to A.S. Mikhailov for helpful discussions.
M. Cloet was funded by KU Leuven grant STG/19/007
L. Arno was funded by a FWO-Flanders fellowship, grant 117702N. D. Kabus is supported by KU Leuven grant GPUL/20/012. H. Dierckx and computational infrastructure were funded by KU Leuven grant STG/19/007.
J. Van der Veken is supported by the Research Foundation–Flanders (FWO) and the National Natural Science Foundation of China (NSFC) under collaboration project G0F2319N, by the KU Leuven Research Fund under project 3E210539 and by the Research Foundation- Flanders (FWO) and the Fonds de la Recherche Scientifique (FNRS) under EOS Projects G0H4518N and G0I2222N.
A. V. Panfilov’s research at Sechenov University was financed by the Ministry of Science and Higher Education of the Russian Federation within the framework of state support for the creation and development of World-Class Research Centers “Digital biodesign and personalized healthcare” (Grant No 075-15-2022-304).

\end{acknowledgments}

%


\bibliography{apssamp}

\end{document}


\appendix

\section{Simulation details and additional visualizations}

The simulation results given in the paper rely on four simulations, of which the details are given here. All spatial quantities are in units of $\si{\milli\meter}$.

\subsection{SIM 1: Basic superfilament}

This simulation initiates a planar superfilament that will be reused in SIM 3 and SIM 4. The simulation domain has physical size $[0,61]^4$, and is discretized with resolution $h=\SI{1}{\milli\meter}$. A superscroll with planar superfilament was initiated by using the crossfield stimulation protocol, with rectangular ${S_1 = [0,6]\times[0,61]\times[0,61]\times[0,61]}$ at $t_1 = \SI{0}{\milli\second}$ and a rectangular ${S_2 = [0,61]\times[0,12]\times[0,61]\times[0,61]}$ at ${t_2 = \SI{15}{\milli\second}}$. The simulation is run with time step $dt=\SI{0.05}{\milli\second}$ using 16 threads.

\subsection{SIM 2: Spherical superfilament}

This simulation corresponds to Fig.~2a from the paper. 
The simulation domain has physical size $[0,61]^4$, and is discretized with resolution $h=\SI{1}{\milli\meter}$. A medium at rest was excited using the crossfield stimulation protocol, with a spherical $S_1$, i.e., the region inside a hypersphere with center $(-20,-20,31,45)$ and radius equal to $40$ at $t_1=\SI{0}{\milli\second}$ and a rectangular ${S_2 = [0,10]\times[0,61]\times[0,61]\times[0,61]}$ at $t_2 = \SI{15}{\milli\second}$. The simulation is run with time step $dt=\SI{0.05}{\milli\second}$ using 16 threads.

Fig.~\ref{fig:supp-mat-2a} shows 2D intersections of the medium at different snapshots of the simulation. See Supplemental Material at [URL will be inserted by publisher] for an animated gif of the same 2D intersections, as well as an interactive version of Fig. 2a from the paper and the corresponding animation.

\subsection{SIM 3: Punctured superfilament}

This simulation corresponds to Fig.~2b from the paper. 
To create a hole in the planar superfilament from SIM 1, a snapshot of part of the medium is taken at $t_3=\SI{40}{\milli\meter}$, i.e., the subdomain ${[0,15]}^4$. In this snapshot, the $u$-value in the region ${S_3 = \{(x,y,z,a)\mid {(\frac{z-8}{5})}^2+{(\frac{a-8}{5})}^2\leq 1\}}$ is set to $u=0$. The spatial resolution is increased to $h=\SI{0.5}{\milli\meter}$. The simulation is run with time step $dt=\SI{0.025}{\milli\second}$ using 16 threads.

Fig.~\ref{fig:supp-mat-2b} shows 2D intersections of the medium at different snapshots of the simulation. See Supplemental Material at [URL will be inserted by publisher] for  an animated gif of the same 2D intersections, as well as an interactive version of Fig. 2b from the paper and the corresponding animation.

\subsection{SIM 4: Anchored superfilament}

This simulation corresponds to Fig.~2c from the paper. 
As initial condition, a snapshot of part of the medium is taken from SIM 1 at $t_3=\SI{40}{\milli\second}$, i.e., the subdomain ${[10,36]\times[2,28]\times[18,34]\times[18,34]}$. Part of the medium is given conductivity zero to make the superfilament anchor to conduction blocks. To construct this non-conducting regions, a 4D hypersphere with radius $\SI{2}{\milli\meter}$ is translated along two lines, given by the equations

\begin{align*}
    l_1&\leftrightarrow \begin{cases} 
        x=6k+10\\ 
        y=-6k+14\\
        z=25k\\
        a=0
        \end{cases},\, 0\leq k \leq 1\,, \\
    l_2&\leftrightarrow \begin{cases}
        x=-6k+16\\
        y=6k+10\\
        z=25k\\
        a=25
    \end{cases},\, 0\leq k \leq 1\,.
\end{align*}

The spatial resolution is increased to $h=\SI{0.5}{\milli\meter}$. The simulation is run with time step $dt=\SI{0.025}{\milli\second}$ using 16 threads.

See Supplemental Material at [URL will be inserted by publisher] for an interactive version of Fig. 2c from the paper.

 \begin{figure*}
     \centering
     \includegraphics[width=0.7\textwidth]{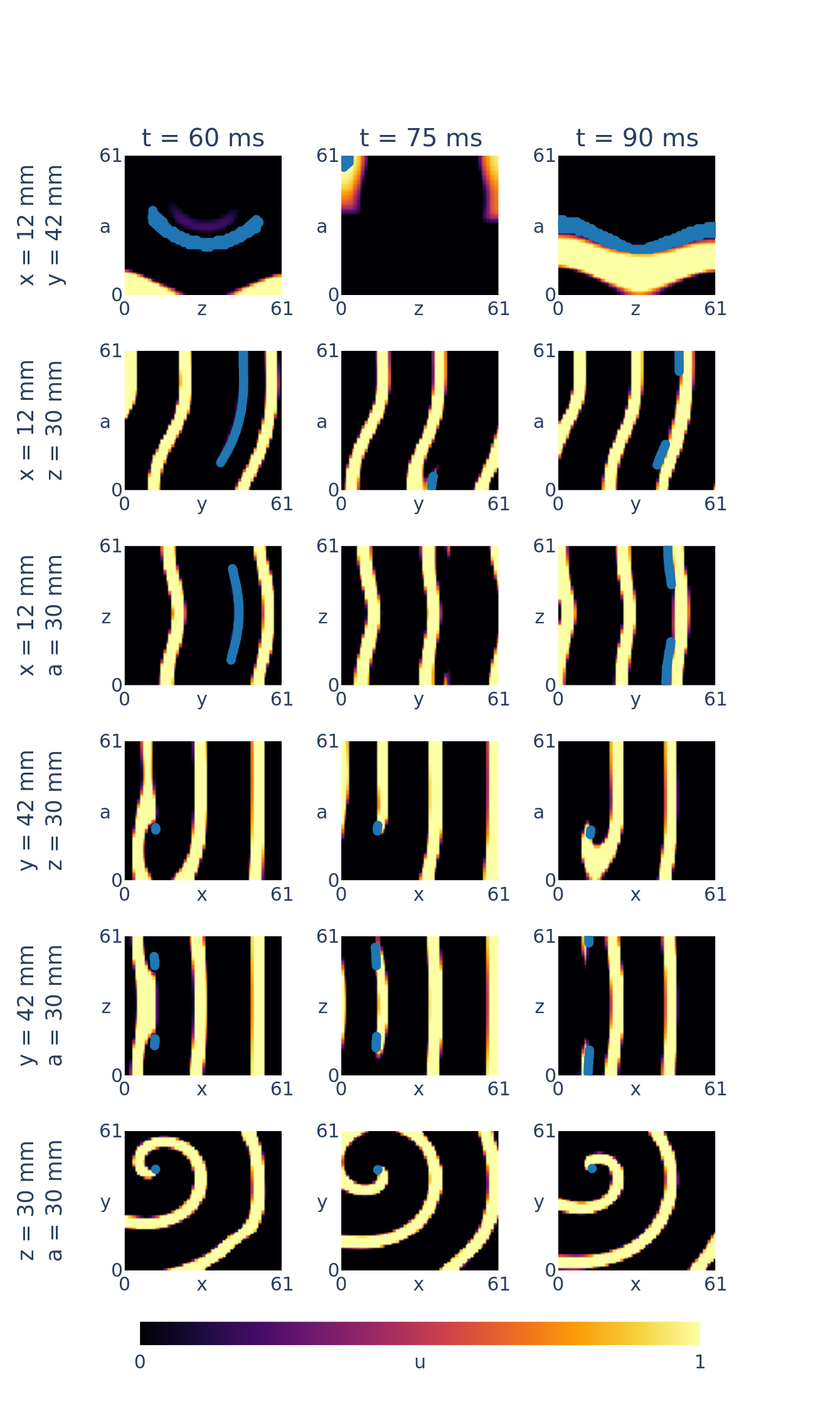}
     \caption{2D intersections of all possible coordinate planes for SIM 2 for three snapshots of the simulation. The coordinates of the intersecting planes are given at the left and the time corresponding to the snapshots at the top. The intersection of the superfilament is given in blue.}
     \label{fig:supp-mat-2a}
 \end{figure*}

 \begin{figure*}
     \centering
     \includegraphics[width=0.7\textwidth]{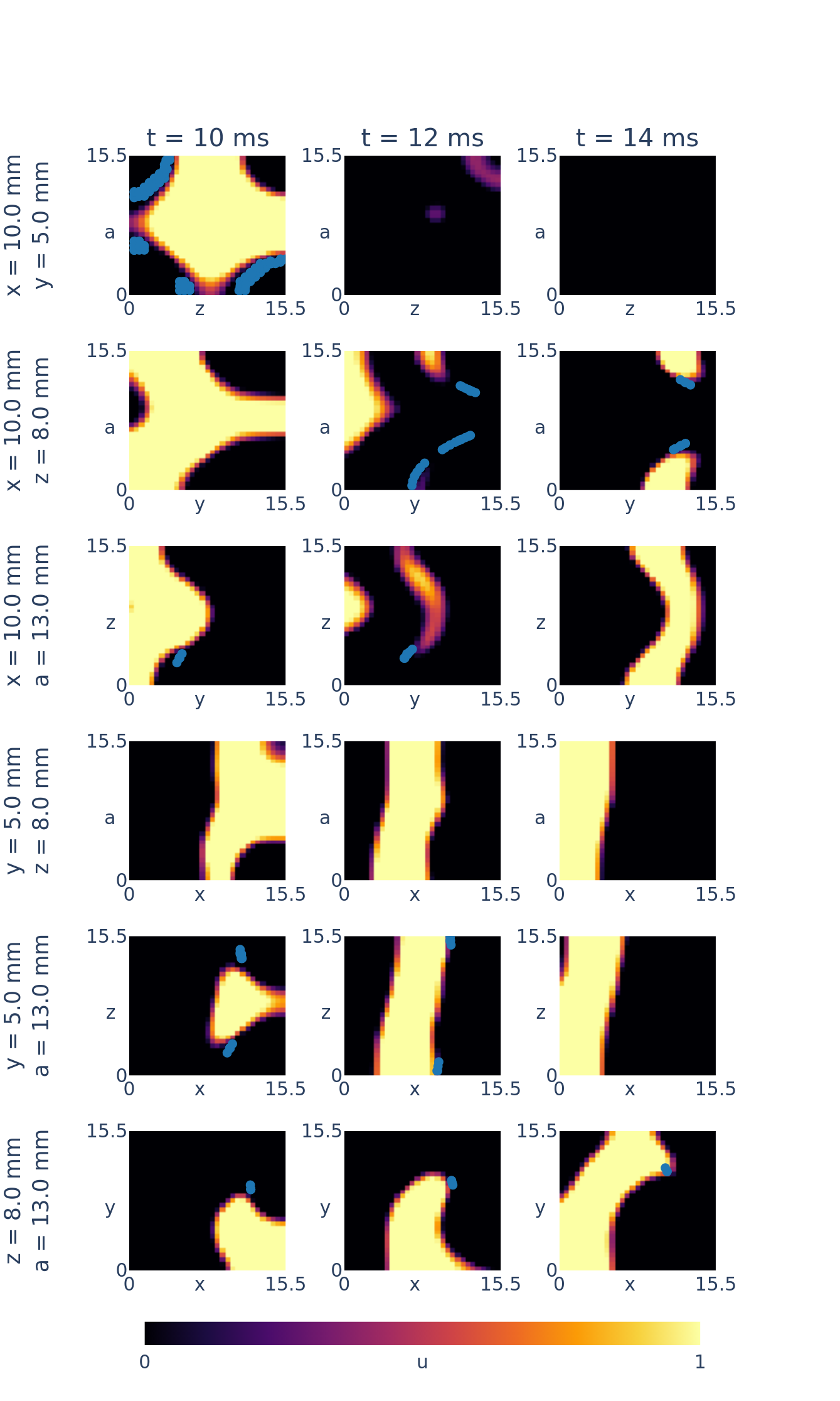}
     \caption{2D intersections of all possible coordinate planes for SIM 3 for three snapshots of the simulation. The coordinates of the intersecting planes are given at the left and the time corresponding to the snapshots at the top. The intersection of the superfilament is given in blue.}
     \label{fig:supp-mat-2b}
 \end{figure*}